\documentclass[12pt]{amsart}
\usepackage{a4}
\usepackage{amssymb}
\usepackage{amsmath}
\usepackage{amsthm}
\usepackage{amstext}
\usepackage{amscd}
\usepackage{latexsym}
\usepackage{graphics}
\usepackage{color}

\swapnumbers
\theoremstyle{plain}
\newtheorem{thm}{Theorem}[section]

\newtheorem{lem}[thm]{Lemma}

\newtheorem{prop}[thm]{Proposition}
\newtheorem{conj}[thm]{Conjecture}
\theoremstyle{definition}
\newtheorem{rem}[thm]{Remark}
\newtheorem{point}[thm]{}

\newcommand{\ilim}{\mathop{\varprojlim}\limits} 
\newcommand{\dlim}{\mathop{\varinjlim}\limits}  

\newcommand{\Intersection}{\bigcap}
\newcommand{\Union}{\bigcup}
\newcommand{\intersection}{\cap}

\newcommand{\Spec}{{\rm Spec \,}}

\renewcommand{\iff}{\mbox{ $\Longleftrightarrow$ }}

\newcommand{\sE}{{\mathcal E}}

\newcommand{\C}{{\mathbb C}}

\newcommand{\F}{{\mathbb F}}

\renewcommand{\P}{{\mathbb P}}

\input{xy}
\xyoption{all}

\begin{document}
\title{Products of Brauer Severi Surfaces}
\author{Amit Hogadi}

\begin{abstract}
Let $\{P_i\}_{1 \leq i \leq r}$ and $\{Q_i\}_{1 \leq i \leq r}$ be two collections of Brauer Severi surfaces (resp. conics) over a field $k$. We show
that the subgroup generated by $P_i's$ in $Br(k)$ is the same as the subgroup 
generated by $Q_i's$ \iff $\Pi P_i $ is birational to $\Pi Q_i$. Moreover in this case 
$\Pi P_i$ and $\Pi Q_i$  represent the same class in $M(k)$, the Grothendieck ring of $k$-varieties. The converse holds if $char(k)=0$. Some of the above implications also hold over a general noetherian base scheme. 
\end{abstract}

\maketitle


\section{Introduction}

\begin{point}[Notation]
Let $S$ denote a noetherian base scheme.  All products, unless otherwise
mentioned, will be over $S$.  The class of any Brauer Severi scheme $P$ over $S$ in $Br(S)$ (the Brauer group of $S$) will be denoted by $P$ itself. 
For a collection of Brauer-Severi schemes $\{P_i\}_{i\in I}$ over $S$, the subgroup generated by the $P_i's$ in $Br(S)$ will be denoted by $<\hspace{-1mm}\{P_i\}_{i\in I}\hspace{-1mm}>$.
$M(S)$ will denote the Grothendieck ring of finite type $S$-schemes (see section(\ref{grothendieckring-section})).
\par All schemes considered will be noetherian.  By a closed subscheme 
we will always mean a reduced closed subscheme. 
\end{point}

\noindent The main result of this paper is the following.
\begin{thm}\label{mainthm} Let $\{P_i\}_{1\leq i \leq r}$ and $\{Q_j\}_{1\leq j\leq r}$ be two collections Brauer Severi surfaces (resp. conics) over $S$.
Consider the following conditions. 
\begin{enumerate}
\item[(i)]$<\hspace{-1mm}\{P_i\}\hspace{-1mm}>=<\hspace{-1mm}\{Q_j\}\hspace{-1mm}>$ in $Br(S)$.
\item[(ii)] $ [\Pi  P_i] = [\Pi_j Q_j]$ in $M(S)$.
\item[(iii)] $\Pi P_i$ and $\Pi Q_j$ are birational.
\end{enumerate}
Then $(i)\Rightarrow (ii)$. If $S$ is reduced then $(i)\Rightarrow (iii)$.
If $S$ is a separated regular scheme then $(i) \iff (iii)$. If $S$ is a separated regular scheme with characteristic zero generic points, 
then $(i)\  (ii) \ \text{and} \  (iii)$ are equivalent. 
\end{thm}
This result has been inspired by \cite{kollar} where 
relations between products of conics in the Grothendieck ring were 
studied for the first time. The above theorem was proved in \cite{kollar}
for conics in the case when $S=\Spec(k)$ where $k$ is a number field or 
function field of an algebraic surface over $\C$. \\
The proof presented here is by induction on $r$.
Working over a general noetherian base scheme $S$ instead of a field   
enables us to run the induction more smoothly.\\

\noindent Recall the following conjecture of Amitsur. 
\begin{conj}[\cite{amitsur}]
Let $k$ be a field and $P$ and $Q$ be $n$-dimensional Brauer-Severi varieties 
over $k$. Then $P$ is birational to $Q$ \iff $P$ and $Q$ generate the same
subgroup in $Br(k)$.
\end{conj}

This conjecture is still unknown in general, however the following 
special cases are known.
\begin{enumerate}
\item  $P$ is split by a cyclic extension  (which is always 
true if $k$ is a local or global field) (see \cite{amitsur}).
\item $index(P)< dim(P)+1$ (see \cite{roquette}).
\item $P=-Q$ in $Br(k)$ (this proves the conjecture for Brauer-Severi surfaces) (see \cite{roquette}).
\item $P= 2Q$ in $Br(k)$ (see \cite{tregub}).
\end{enumerate}
\begin{rem}
In addition to Brauer-Severi surfaces and conics, the proof of (\ref{mainthm}) 
presented here also works for Brauer-Severi varieties of prime index if one 
assumes Amitsur's conjecture for this case.
\end{rem}
\section{Preliminaries on the Grothendieck Ring}\label{grothendieckring-section}
\begin{point}(Grothendieck Ring).
Let $S$ be any scheme. Let $M(S)$ denote the free abelian group 
generated on isomorphism classes of reduced finite type $S$-schemes modulo 
the relations 
$$ [X]=[U]+[Z] $$
where $X$ is a reduced $S$ scheme and $U\subset X$ is an open subset with complement 
$Z$ (with reduced scheme structure). For any $S$ scheme $X$, we will use the 
notations $[X]_S$ or just $[X]$ to denote the class of $X^{red}$ in $M(S)$. 
For $S$ schemes $X,Y$ define 
$$[X]_S\cdot [Y]_S = [(X\times_SY)]_S$$
This makes $M(S)$ into a commutative and associative ring with $[S]$ being 
the identity in this ring. $M(S)$ is called 
the Grothendieck ring of finite type $S$-schemes. Notice that $M(S)$ depends only on the reduced structure of $S$.
\end{point}
\begin{point}($f^*$ and $f_*$). Given any morphism $f:T\to S$, there is functorial ring homomorphism
$f^*:M(S)\to M(T)$ induced by base extension $X\to X\times_ST$.
Moreover if $f$ is itself of finite type, one also has a morphism of $M(S)$-modules
$f_*:M(T)\to M(S)$ induced by considering any $T$-scheme as an $S$-scheme via $f$.\\

\par Suppose we have a filtered inverse system of schemes $\{S_i\}_{i\in I}$ 
such that the inverse limit $\ilim S_i$ exists.
Then one gets a natural ring homomorphism  
$$ \dlim M(S_i) \to M(\ilim S_i)$$
The following special case is of special interest.
\end{point}

\begin{prop}\label{prop:mlimit}
Let $S$ be an integral scheme. Let $\{ U\}_{U\subset S}$ be the (filtered) 
inverse system of nonempty
open sets of $S$. Let $K$ be the function field of $S$. Then the natural ring homomorphism
$ \mathop{\dlim}_{U\subset S} M(U) \to M(\Spec(K))$ is an isomorphism.
\end{prop}
\begin{proof}
Any finite type $K$-scheme $X_K$, is the generic fibre of some finite type 
$U$-scheme $X_U$ for some nonempty open set $U$ of $S$. This shows the above map 
is surjective. After shrinking $U$ if necessary, 
any closed subscheme $Z_K\subset X_K$ can be realized as the generic fibre 
of a closed subscheme $Z_U \subset X_U$. This shows the map is injective.
\end{proof}
\noindent Finally, we recall the following elementary proposition.
\begin{prop}[well-known]\label{prop:vb}
Let $\sE$ be a vector bundle on $S$ of rank $n+1$. Then 
$[\P roj(\sE)]= [\P^n_S]$ in $M(S)$.
\end{prop}
\begin{proof}
Let $S=\Union_{i=1}^m U_i$ be an open cover such that $\sE_{|U_i}$ is trivial for each $i$. We now proceed by induction on $m$. If $m=1$, $\sE$ is a trivial vector bundle and the statement is obvious. For $m>1$, let $S'=\Union_{i=1}^{m-1} U_i$. Then by induction $[\P roj(\sE_{|S'})] =[\P^n_{S'}]$. Since $\sE$ is trivial on $U_m$ and $S'\Intersection U_m$, we also have $[\P roj(\sE_{|U_m})]=[\P^n_{U_m}]$ and $[\P roj(\sE_{|S'\Intersection U_m})]=[\P^n_{S'\Intersection U_m}]$. The result now follows from the following equalities in $M(S)$. 
\begin{align*}
 [\P roj(\sE)]  & = [\P roj(\sE_{|S'})] + [\P roj(\sE_{|U_m})] - [\P roj(\sE_{|S'\Intersection U_m})] \\
 [\P^n_S]  & = [\P^n_{S'}] + [\P^n_{U_m}] - [\P^n_{S'\Intersection U_m}] 
\end{align*}
\end{proof}

\section{The Cremona Map}
\begin{point}(The Cremona Map) \label{p:cremona} Let $K$ be a field. 
Let us recall the following well known birational map from 
$\P^2_K$ to itself. 
$$ \phi:\P^2_K \dashrightarrow \P^2_K \hspace{4mm} [X,Y,Z]\to [YZ,XZ,XY]$$
$\phi$ can be defined everywhere on $\P^2_K$ outside the reduced closed subscheme
$$ B= \{[1,0,0],[0,1,0], [0,0,1]\}$$
Let $X\stackrel{p}{\to}\P^2_K$ be the blow up of $\P^2_K$ with center $B$. Then $\phi$ defines a morphism 
 $X\stackrel{q}{\to} \P^2_K$ such that the following diagram commutes.
$$
\xymatrix{
        & X \ar[ld]_p\ar[rd]^q & \\
\P^2_K \ar@{-->}[rr]^{\phi} &      & \P^2_K 
}
$$
One can check that $q:X\to \P^2_K$ is itself is the blowup of $\P^2_K$
again with center $B$.
\end{point}

\noindent The following result was essentially proved in \cite{roquette}.
\begin{thm}\label{cremona1}
Let $P$ and $Q$ be Brauer Severi surfaces over a field $K$.
Assume that $P=2Q$ in $Br(K)$. Then there exists a birational 
map $\phi:P\dashrightarrow Q$ which after going to $\overline{K}$
(the separable closure of $K$) is isomorphic to the Cremona map.
\end{thm}

\begin{thm}\label{cremona2}
Let $K$ be any field and let $P$ and $Q$ be Brauer-Severi surfaces which 
generate the same subgroup in $Br(K)$. Then $[P]=[Q]$ in $M(K)$.
\end{thm}
\begin{proof}
Let $ \phi: P \dashrightarrow Q $ be a map as guaranteed by Theorem(\ref{cremona1}).
Let $B$ (resp. $B'$) be the base locus of the map $\phi$ (resp. $\phi^{-1}$). 
Without loss of generality we may assume that $P$ defines a nontrivial class 
in $Br(K)$ and thus has no $K$-point. Then $B$ (resp. $B'$) is a closed $L$-point (resp. $L'$-point) of $P$ (resp. $Q$) for some degree $3$ separable field extension $L/K$ (resp. $L'/K$). We claim that $L/K$ and $L'/K$ are isomorphic field extensions. Let $X$ be the blow up
of $P$ at $B$. Then we have the following Hironaka hut.
$$
\xymatrix{
        & X \ar[ld]_p\ar[rd]^q & \\
P\ar@{-->}[rr]^{\phi} &      & Q 
}
$$
Here $q$ is the blowup of $Q$ at $B'$ (see (\ref{p:cremona})).
To show that $L/K$ and $L'/K$ are isomorphic it is enough to show that $B'\times_K \Spec(L)$ has an $L$-point. But after base extending to $L$, $P_L=P\times_KL \cong \P^2_L$ and $p$ is the blow up of three $L$-points of $P_L$, say $x,y,z$. Let $L_{xy}$ be the unique 
line in $P_L$ joining $x$ and $y$ and similarly for $L_{yz},L_{xz}$. Then 
the birational transform of $L_{xy}\Union L_{yz}\Union L_{xz}$ is the exceptional
locus of $q$. Thus the image, $B'\times_K\Spec(L)$, of this  exceptional locus is the disjoint union of $3$ points. This proves the claim.
\par Thus as $K$-varieties, the exceptional locus of $p$ (resp. $q$) is 
isomorphic to $\P^1_L$. Thus $ [P]= [X]-[\P^1_L]+[\Spec(L)] = [Q]$.
\end{proof}

\section{Proof of the main theorem}
\noindent For any morphism $f:T\to S$ and any $S$-scheme $X$ let $X_T=X\times_ST$.
\begin{lem}\label{cor:dim}
Let $X$ be any finite type $S$-scheme. Let $P$ be a Brauer-Severi scheme over $S$
of relative dimension $n$.
Assume that the class of $P$ in $Br(S)$ lies in the kernel of $Br(S)\to Br(X)$. 
Then $$[X\times P]_S= [\P^n_X]_S$$
\end{lem}
\begin{proof}
It is enough to prove $[X\times P]_X= [\P^n_X]_X$ since the required 
equality can then by obtained by using the natural map $M(X)\to M(S)$.
But by assumption, the class represented by $X\times P$ in $Br(X)$ is zero.
Hence there exists a vector bundle $\sE$ on $X$ such that $X\times P$ is 
isomorphic to $\P roj(\sE)$ as $X$-schemes. The result now follows
from (\ref{prop:vb}).
\end{proof}

\begin{proof}[Proof of (\ref{mainthm})] In order to avoid unnecessary repetition, we will only prove the theorem for Brauer-Severi surfaces. We proceed by induction on $r$.\\
\noindent $(i)\Rightarrow (ii)$: \\
\noindent {\it Step(1)}: One can quickly reduce to proving the statement 
in the case when $S$ is integral.
Let $S=S_1\Union S_2$ be the decomposition of $S$ into 
two closed subschemes. Then for any $S$-scheme $X$ of finite type, we have
$$ [X]_S= [X_{S_1}]_S+[X_{S_2}]_S-[X_{S_{12}}]_S\ \ \ \text{where} \ \ S_{12}=S_1\intersection S_2$$
Hence by noetherian induction and the above formula, 
it is enough to prove the theorem in the case when $S$ is irreducible. 
Moreover the natural ring homomorphism  $M(S)\to M(S^{red})$ is an isomorphism. 
Thus we may assume $S$ is integral.\\

\noindent {\it Step(2)}:($r$=1)
Let $P/S$ and $Q/S$ be two Brauer Severi surfaces. Let $K$ be the function field of $S$. 
Then by (\ref{cremona2}) 
$$[P_K]= [Q_K] \hspace{5mm} \text{in} \hspace{4mm} M(K)$$
By (\ref{prop:mlimit}), there exists a nonempty open set $U$ of $S$ such that 
$$[P_U]= [Q_U] \hspace{5mm} \text{in} \hspace{4mm} M(U)$$ Let $Z$ be the complement of $U$. By noetherian 
induction $[P_Z]=[Q_Z]$. Since $[P]=[P_U]+[P_Z]$
and similarly for $[Q]$, we get that $[P]_S=[Q]_S$.\\

\noindent {\it Step(3)}:
Suppose the dimension of $<\hspace{-1mm}\{P_i\}\hspace{-1mm}>$ (and hence also of $<\hspace{-1mm}\{Q_j\}\hspace{-1mm}>$) as an $\F_3$ vector space is strictly less than $r$. Then without loss of generality we may assume that the class of $P_r$ is contained in the subgroup generated by $\{P_i\}_{1 \leq i \leq r-1}$. Then by Corollary(\ref{cor:dim})
$$ [\mathop{\Pi}_i \ P_i]= [\mathop{\Pi}_{i\leq r-1}P_i \times \P^2_S]$$
And similarly for the $Q_i's$. By induction 
$$[\Pi_{1\leq i\leq r-1} P_i]=[\Pi_{1\leq j \leq r-1} Q_j]$$ which implies 
$$ [\Pi_{1\leq i \leq r} P_i] = [\Pi_{1\leq i\leq r-1} P_i\times \P^2_S]=[\Pi_{1\leq j \leq r-1} Q_j\times \P^2_S] = [\Pi_{1\leq j \leq r} Q_j]$$
Hence it is enough to prove the theorem under the extra assumption that dimension of $<\hspace{-1mm}\{P_i\}_{1\leq i \leq r}\hspace{-1mm}>$ 
as an $\F_3$ vector space is $r$.\\

\noindent {\it Step(4)}: Since class of $Q_1$ is the in the subgroup 
generated by $P_i's$ in $Br(S)$, we have the following equation in $Br(S)$ 
$$Q_1 = \sum a_i P_i \hspace{4mm}, a_i \in \F_3$$
By $Step(3)$, at least one of the $a_i's$ is nonzero. Without loss of generality we may assume $a_1\neq 0$. Thus $Q_1= a_1P_1 + \sum_{i\geq 2} a_i P_i$ in $Br(S)$. We first claim that 
$$[Q_1 \times \mathop{\Pi}_{2\leq i \leq r}P_i ]_S = [P_1 \times \mathop{\Pi}_{2\leq i \leq r}P_i ]_S$$
Put $Y=\mathop{\Pi}_{2\leq i \leq r}P_i$. Now $Q_1\times Y$ and $P_1\times Y$
generate the same subgroup in $Br(Y)$. Hence by $Step(2)$, $[Q_1\times Y]_Y = [P_1\times Y]_Y$.
The claim now follows by using the push forward map $M(Y)\to M(S)$.
Now to prove the theorem it is enough to show 
$$[Q_1 \times \mathop{\Pi}_{2\leq i \leq r}P_i ]_S=[Q_1 \times \mathop{\Pi}_{2\leq i \leq r}Q_i ]_S$$
But the subgroup generated by $\{P_i\}_{2\leq i \leq r}$ in $Br(Q_1)$ is the same as the subgroup 
generated by $\{Q_i\}_{2\leq i \leq r}$. Hence by induction on $r$ we have 
$$ [Q_1 \times \mathop{\Pi}_{2\leq i \leq r}P_i ]_{Q_1}=[Q_1 \times \mathop{\Pi}_{2\leq i \leq r}Q_i ]_{Q_1}$$
Again, the claim follows by using the push forward map $M(Q_1)\to M(S)$.\\

\noindent $(i)\Rightarrow (iii)$: ($S$ is reduced): As in the proof of 
$(i)\Rightarrow (ii)$, we proceed by induction on $r$. 
One first reduces the proof to the case when $S$ is integral. Then by noetherian
induction and the known result for the case when $S$ is the spectrum of a field
we prove case $r=1$. After a possible re-indexing, one then proves $P_1 \times \mathop{\Pi}_{2\leq i \leq r} P_i$
is birational to $Q_1 \times \mathop{\Pi}_{2\leq i \leq r}Q_i$ by first 
comparing $P_1 \times \mathop{\Pi}_{2\leq i \leq r} P_i$
and $Q_1 \times \mathop{\Pi}_{2\leq i \leq r} P_i$ and then comparing 
$Q_1 \times \mathop{\Pi}_{2\leq i \leq r} P_i$ and 
$Q_1 \times \mathop{\Pi}_{2\leq i \leq r} Q_i$. Since the argument is very similar
to the one above, we leave the details to the reader. \\

\noindent ($S$ is a separated regular scheme): Without loss of generality we may assume $S$ is connected. Let $K$ be the function field of $S$. Since $S$ is regular, $Br(S)\to Br(K)$ is injective. Thus in order to prove $(iii) \Rightarrow (i)$ and $(ii)\Rightarrow (i)$ we may replace $S$ by $\Spec(K)$. The kernel of $Br(K) \to Br(\Pi P_i)$ is equal to $<\{P_i\}>$.  Moreover this kernel depends only on the stable birational class of $\Pi P_i$. This proves that $(iii)\Rightarrow (i)$.
\par  If $K$ is of characteristic zero then $(ii)\Rightarrow (i)$ follows from the fact that any two smooth projective varieties having the same image in $M(K)$ are stably birational (see \cite{ll}). 
\end{proof}

\noindent {\bf Acknowledgement}:
I thank my advisor J\'anos Koll\'ar for his encouragement and useful discussions.
I also thank Chenyang Xu for useful comments and discussions.

\bibliographystyle{amsplain}

\begin{thebibliography}{0000000000}

\bibitem{amitsur}
Amitsur S. A., Generic splitting fields of central simple algebras, 
{\it Ann. of Math (2)} {\bf 62}, 1955, 8--43.

\bibitem{kollar}
Koll\'ar J., Conics in Grothendieck Ring , 
 {\it Adv. Math.}  {\bf 198}  (2005), {\bf no. 1}, 27--35.

\bibitem{ll} Larsen, M.; Lunts, V.A.; Motivic measures and stable birational geometry.  
{\it Mosc. Math. J.}  {\bf 3}  (2003),  {\bf no. 1}, 85--95. 

\bibitem{roquette}
Roquette P., On the Galois cohomology of the projective linear group and its
applications to the construction of generic splitting fields of algebras, 
{\it Math. Ann.} {\bf 150}, 1963, 411-439.

\bibitem{tregub}
Tregub, S. L., Birational equivalence of Brauer-Severi manifolds,
{\it Uspekhi Mat. Nauk} {\bf 46}, 1991, 217--218; English translation in 
{\it Russian Math. Surveys} {\bf 46}, 1992, 229. 


\end{thebibliography}

\vspace{7mm}
\noindent  email : amit@math.princeton.edu

\end{document}